\documentstyle[12pt]{article}     
\input amssym.def
\input amssym.tex

\def\({{\em (}}
\def\){{\em )}}
\def\:{{\em\,:}}

\def\Section#1{\vspace{30truept}\addtocounter{section}{1}\centerline{\Large\bf
	\arabic{section}.~~#1}\vspace{12pt}}
\newtheorem{definition}{Definition}
\newtheorem{theorem}{Theorem}
\newtheorem{lemma}{Lemma}
\newtheorem{proposition}{Proposition}
\newtheorem{corollary}{Corollary}

\newtheorem{remark}{Remark}

\title{L-CONNECTIONS AND ASSOCIATED TENSORS}
\author{Nabil L. Youssef} 
\date{}

\begin{document}                     
\bibliographystyle{plain}
\maketitle                           
\vspace{-1.3cm}
\begin{center}
{Department of Mathematics, Faculty of Science,\\ Cairo
University, Giza, Egypt.} 
\end{center}
\vspace{1cm}
\maketitle                         
\smallskip
\Section {Introduction}
\par The theory of connections in Finsler geometry is not satisfactorily 
established as in Riemannian geometry. Many trials have been carried 
out to build up an adequate theory. One of the most important  
in this direction is that of Grifone ([3] and [4]). His approach to the
theory of nonlinear connections was accomplished in [3], in which 
his new definition of a nonlinear connection is easly handled from the
algebraic point of view. Grifone's approach is based essentially
on the natural almost-tangent structure $J$ on the tangent
bundle $T(M)$ of a differentiable manifold $M$. This structure was
introduced and investigated by Klein and Voutier [5]. Anona in
[1] generalized the natural almost-tangent
structure by considering a vector $1$-form $L$ on a manifold $M$
(not on $T(M)$) satisfying certain conditions. He investigated 
the $d_L$-cohomology induced on $M$ by $L$ and
generalized some of Grifone's results.
\par In this paper, we adopt the point of view of Anona [1] to
generalize Grifone's theory of nonlinear connections [3]: We
consider a vector $1$-form $L$ on $M$ of constant rank such that
$[L,L]=0$ and that $Im(L_z)=Ker(L_z)$; $z\in M$. We found
that $L$ has properties similar to those of $J$, which enables
us to generalize systematically the most important results of
Grifone's theory.
\par The theory of Grifone is retrieved, as a special case of
our work, by letting $M$ be the tangent bundle of a
differentiable manifold and $L$ the natural almost-tangent
structure $J$. 

\Section {Notations and Preliminaris}
\par The following notations will be used throughout the paper:\newline
$M$: a differentiable manifold of class $C^\infty$ and of finite dimension.\newline
$T(M)$: the tangent bundle of $M$.\newline
${\goth F}(M)$: the $\Bbb R$-algebra of differentiable functions on $M$.\newline
${\goth X}(M)$: the ${\goth F}(M)$-module of vector fields on $M$.\newline
$\Phi^p(M)$: the ${\goth F}(M)$-module of scaler $p$-forms on $M$ ($\Phi^0(M)={\goth F}(M))$.\newline
$\Psi^l(M)$: the ${\goth F}(M)$-module of vector $l$-forms on $M$ ($\Psi^0(M)={\goth X}(M))$.\newline
$\Phi(M)$, $\Psi(M)$: the corresponding graded rings.\newline
$J$: the natural almost-tangent structure on $T(M)$ ([5] and [3]).\newline
${\goth L}_X$: the Lie derivative with respect to $X\in{\goth X}(M).$\newline
$i_K$: the interior product with respect to $K\in\Psi(M)$.
\par All geometric objects considered in this paper are supposed to
be of class $C^\infty$. The formalism of Fr\"{o}licher-Nijenhuis [2]
will be the fundamental tool of the whole work.

Let $M$ be a paracompact differentiable manifold of dimension $m+n$. 
Let\linebreak $L\in\Psi^1(M)$ be a vector $1$-form on $M$ of constant rank
$m\ge1$ and whose \linebreak Nijenhuis tensor $N_L=\frac 12
[L,L]$ is zero. The vector form $L$ defines on $M$ a distribution ${\goth D}:
z\in M\longmapsto L_z(T_z(M))$ which is completely integrable. The
manifold $M$ equipped with this structure is said to be foliated
by $L$ (vari\'{e}t\'{e} feuillet\'{e}e par $L$) [6].
\par For any chart of $M$ with domain $U$ we have a system of
local coordinates $(x^1,\ldots, x^n; y^1,\ldots, y^m)$ such that
$(\frac{\partial}{\partial y^1},\ldots, \frac{\partial}{\partial y^m})$
is a basis for $LT(U))$. Let $\alpha, \beta,\ldots\in\{1,\ldots,
n\}$ and $i, j,\ldots\in\{1,\ldots, m\}$, then the distribution $\goth D$
is defined by the \linebreak equations
\begin{equation}
dx^\alpha=0
\end{equation}
We can define a supplementary distribution ${\goth D}^s$ by the equations 
\begin{equation}
\theta^i=dy^i+\Gamma^i_\alpha dx^\alpha=0
\end{equation}
where $\Gamma^i_\alpha$ are $C^\infty$ functions. We thus have
the decomposition
\begin{equation}
T_z(M)={\goth D}_z\oplus {\goth D}^s_z\qquad \forall z\in M
\end{equation}
If $v:T_z(M)\longrightarrow{\goth D}_z$ and $h:T_z(M)\longrightarrow{\goth D}^s_z$
are the associated projectors, the vector $1$-form $\Gamma=h-v$
defines an almost-product structure on $M$. The curvature $R$ of
$\Gamma$ being defined by $R=-\frac 1{2}[h,h]$, it can be proved
that the distribution ${\goth D}^s$ is completely integrable if,
and only if, $R$ vanishes.
\par Let $C$ be the vector field on $M$ generating the
one-parameter group $h_t$ on $U$, 
$$h_t(x^\alpha,y^i)=(x^\alpha,e^ty^i);\quad t\in\Bbb R.$$
Hence, $C$ can be expressed on $U$ in the form:
\begin{equation}
C=y^i\frac{\partial}{\partial y^i}
\end{equation}
It is clear that $C$ can be defined on the whole manifold $M$.
The vector field $C$ is called the canonical vector field on $M$.\newline
\par The following definitions will be used in the sequel ([1] and [3]):
\begin{definition}\quad
A scalar form $\omega\in\Phi(M)$ is said to be homogeneous of
degree r if ${\goth L}_C\omega=r\omega$.
\par A vector form $K\in\Psi(M)$ is said to be homogeneous of
degree r if $[C,K]=(r-1)K$.
\end{definition}
\begin{definition}\quad
A scalar form $\omega\in\Phi(M)$ is said to be $L$-semibasic
if $\,i_X\omega=0\,$ for all vector field $X\in Ker(L)$.
\par A vector form $K\in\Psi(M)$ is said to be $L$-semibasic
if $\,LK=0\,$ and $\,i_XK=0\,$ for all vector field $X\in Ker(L)$.
\end{definition}
\begin{definition}\quad
A vector field $S\in{\goth X}(M)$ is said to be an $L$-semispray
on $M$ if $LS=C$.
\par An $L$-semispray $S$ on $M$ is said to be an $L$-spray on
$M$ if $S$ is homogeneous of degree $2$: {\em(}$[C,S]=S${\em)}.
\end{definition}
\begin{definition}\quad
Let $K$ be an $L$-semibasic scalar or vector $k$-form ($k\ge1$)
on $M$. The potential $K^\circ$ of $K$ is the $L$-semibasic
$(k-1)$-form defined by $\,K^\circ=i_SK\,$, where $S$ is an
arbitrary $L$-semispray. 
\par Clearly, $K^\circ$ does not depend on the choice of the $L$-semispray $S$.
\end{definition}

\par Finally, we have [1] for all $X\in{\goth X}(M)$,
\begin{equation}
L[LX,S]=LX
\end{equation}
\begin{equation}
L[S,S']=S-S'
\end{equation}
where $S$ and $S'$ are arbitrary $L$-semisprays.

\Section {$L$-Connections}
In addition to the conditions imposed on $L$ in section 1, we
assume further that the image of $L$ coincides with its kernel.
It follows immediately that:\newline
--- $m=n$: the dimension of $M$ is even and equals $2n$.\newline
--- $L^2=0$: $L$ is an almost-tangent structure on $M$.\newline
--- ${\goth D}_z=Im(L_z)=Ker(L_z)$.
\par Moreover, we have
\begin{lemma}\quad
The vector form $L$ is homogeneous of degree zero{\em:} $[C,L]=-L.$
\end{lemma}
{\bf Proof.}
Let $(U,\phi)$ be a chart of $M^{2n}$ with natural local coordinates 
$(x^i,y^i)$;\linebreak $i=1,\ldots, n$. Any vector field $X\in{\goth X}(M^{2n})$
is then expressed on $U$ by:
$$X= X^i\frac{\partial}{\partial x^i}+ 
Y^i\frac{\partial}{\partial y^i}$$ 
As $L(\frac{\partial}{\partial x^i})=\frac{\partial}{\partial y^i}$
and $L(\frac{\partial}{\partial y^i})=0$, then 
$$LX=X^i\frac{\partial}{\partial y^i}$$
The canonical vector field $C$ being given by (4), we have
$$[C,LX]=y^i\frac{\partial X^j}{\partial y^i}\frac{\partial}{\partial y^j}- X^j\frac{\partial}{\partial y^j},\quad
L[C,X]=y^i\frac{\partial X^j}{\partial y^i}\frac{\partial}{\partial y^j}$$
The above two equalities imply that
$$[C,L]X=[C,LX]-L[C,X]=-X^j\frac{\partial}{\partial y^j}=-LX $$
Hence the result.\quad $\Box$
\begin{lemma}{\em [3]}\quad
Let $K$ be an $L$-semibasic vector $k$-form, homogeneous of
degree $r$ with $r+k\not=0$. Then we have
$$ K=\frac 1{r+k}(\,[L,K]^\circ+[L,K^\circ]\,)$$
\end{lemma}
\begin{definition}\quad
A vector $1$-form $\Gamma$ on $M^{2n}$ is called an $L$-connection
on $M^{2n}$ if\linebreak $L\Gamma=L$ and $\Gamma L=-L$.
\par An $L$-connection $\Gamma$ on $M$ is said to be homogeneous
if it is homogeneous of degree $1$ as a vector form{\em:} $[C,\Gamma]=0$.
\end{definition}
\begin{theorem}\quad
A vector $1$-form $\Gamma$ on $M^{2n}$ is an $L$-connection if,
and only if, $\Gamma$ defines an almost-product structure on
$M^{2n}$ such that, for all $z\in M^{2n}$, the eigenspace of
$\Gamma_z$ corresponding to the eigenvalue $(-1)$ coincides with
the image space of $L_z$.
\end{theorem}
{\bf Proof.}\newline
Necessity. Let $\Gamma\in\Psi^1(M)$ be an $L$-connection on $M$.
We have:\newline
$L\Gamma=L\iff L(\Gamma-I)=0\iff Im(\Gamma-I)\subset Ker(L)=Im(L),$\newline
$\Gamma L=-L\iff (\Gamma+I)L=0\iff Im(L)\subset Ker(\Gamma+I).$\newline
It follows that $\,Im(\Gamma-I)\subset Ker(\Gamma+I)$. Hence,
$(\Gamma-I)(\Gamma+I)=0$ or $\Gamma^2=I$.
\par Now, let $V_{\Gamma_z}(-1)$ be the eigenspace of $\Gamma_z$
corresponding to the eigenvalue $(-1)$. Let $X\in T_z(M)$, then
(we drop the suffix $z$ for simplicity):\newline
$X\in V_\Gamma(-1)\Longrightarrow \Gamma X=-X\Longrightarrow
LX=L(\Gamma X)=-LX\Longrightarrow  X\in Ker(L)=Im(L)$,
$X\in Im(L)\Longrightarrow X=LY;\; Y\in T_z(M)\Longrightarrow 
\Gamma X=\Gamma(LY)=-X\Longrightarrow 
X\in V_\Gamma(-1).$\newline
Hence, $V_{\Gamma_z}(-1)=Im(L_z).$\newline
Sufficiency. Let $\Gamma\in\Psi^1(M)$ be such that $\Gamma^2=I$ and
$V_{\Gamma_z}(-1)=Im(L_z) \quad\forall z\in M.$\linebreak
$\Gamma^2=I\Longrightarrow(\Gamma+I)(\Gamma-I)=0\Longrightarrow 
Im(\Gamma-I)\subset Ker(\Gamma+I)=Im(L)=Ker(L)\Longrightarrow
L(\Gamma-I)=0\Longrightarrow L\Gamma=L$. 
On the other hand, for all $X\in T_z(M)$,\newline
$LX\in Im(L)=V_\Gamma(-1)\Longrightarrow
\Gamma(LX)=-LX\Longrightarrow \Gamma L=-L.$\newline
Hence, $\Gamma$ is an $L$-connection on $M$.\quad $\Box$
\medskip
\par We set:\newline
$h=\frac12(I+\Gamma)$: horizontal projector associated with
$\Gamma$ ($h^2=h$).\newline
$v=\frac12(I-\Gamma)$: vertical projector associated with
$\Gamma$ ($v^2=v$).\newline
Then, the decomposition (3) becomes:
\begin{equation}
T_z(M)=V_z(M)\oplus H_z(M)\qquad \forall z\in M 
\end{equation}
where $V_z(M)={\goth D}_z=V_{\Gamma_z}(-1)=Im(L_z)=Im(v_z)=Ker(h_z)$:
vertical space at $z$,\linebreak
\phantom{mmm} $H_z(M)={\goth D}^s_z=V_{\Gamma_z}(+1)=Im(h_z)=Ker(v_z)$:
horizontal space at $z$.\newline 
And the functions $\Gamma^i_\alpha$ in (2) are the coefficients of
the $L$-connection $\Gamma$.
\par In terms of vector bundles, (7) takes the form: 
$$T(M)=V(M)\oplus H(M) $$
\par The projectors $h$ and $v$ have the properties:
\begin{equation}
Lh=L,\quad hL=0,\quad \Gamma h=h\Gamma=h;\quad 
Lv=0,\quad vL=L,\quad \Gamma v=v\Gamma=-v 
\end{equation}
\par Moreover, we have from (5) and (8):
\begin{equation}
h[LX,S]=hX
\end{equation}
\begin{remark}\quad
{\em Definition 5 generalizes the connection of Grifone [3].
In fact, if $M^{2n}=T(N)$, where $N$ is a differentiable
manifold of dimension $n$, and if $L=J$; the natural
almost-tangent structure on $T(N)$, then the $L$-connection just
defined is nothing but the connection of Grifone and $V(TN)$ and
$H(TN)$ are the usual vertical and horizontal
bundles respectively [3].}
\end{remark}
\begin{proposition}\quad
To each $L$-connection $\Gamma$ on $M$ there is associated a unique\linebreak
$L$-semispray $S$ such that 
$$\frac12[C,\Gamma]^\circ+S-[C,S]=0.$$
Moreover, if $\Gamma$ is homogeneous, $S$ is an $L$-spray.
\par $S$ is called the canonical $L$-semispray associated with $\Gamma$. 
\end{proposition}
{\bf Proof.}
It suffices to write
\begin{equation}
S=hS',
\end{equation}
where $h$ is the horizontal projector associated with $\Gamma$
and $S'$ is an arbitrary \linebreak $L$-semispray. It is easy to verify
that $S$ has the required properties.
\begin{theorem}~\newline
\(a\) If $S$ is an $L$-semispray on $M$, then $[L,S]$ is an
$L$-connection on $M$ whose canonical $L$-semispray is $\,\frac12(S+[C,S])$.\newline
\(b\) If $S$ is an $L$-spray on $M$, then $[L,S]$ is a homogeneous
$L$-connection on $M$ whose canonical $L$-spray is $S$ itself.
\end{theorem}
{\bf Proof.}\newline
(a) Using (5), taking the properties of $L$ into account, one
can show that $L[L,S]=L$ and $[L,S]L=-L.$ This means that $[L,S]$
is an $L$-connection on $M$.
\par Let $S'$ be an arbitrary $L$-semispray on $M$. Then, using
(6), we get $hS'=\frac12(I+[L,S])S'=\frac12(S'+[LS',S]-L[S',S])
=\frac12(S'+[C,S]-(S'-S))=\frac12(S+[C,S]).$\newline
(b) If $S$ is an $L$-spray on $M$, then $[C,S]=S$. So, applying
the Jacobi identity, we get $[C,[L,S]]=0$. It follows that $[L,S]$
is homogeneous. Moreover, $hS'=S$ for all $L$-semispray $S'$ on
$M$.\quad $\Box$ 

\Section {Torsion of an $L$-Connection}
\begin{definition}\quad
Let $\Gamma$ be an $L$-connection on $M$. The torsion of $\Gamma$
is the vector $2$-form
$$T=\frac12[L,\Gamma]$$
\end{definition}
\par For all $X,Y\in{\goth X}(M)$, we have after some calculations: 
\begin{equation}
T(X,Y)=v[LX,hY]+v[hX,LY]-L[hX,hY]
\end{equation}
\begin{lemma}\quad
The torsion $T$ of an $L$-connection $\Gamma$ on $M$ has the
following properties\:\newline
\(a\) $T$ is $L$-semibasic.\newline
\(b\) $T(hX,hY)=T(X,Y)$.\newline
\(c\) $T=[L,h]=-[L,v]$.\newline
\(d\) $[L,T]=0${\em:} $T$ is $L$-closed.\newline
\(e\) $[C,T]=-T+\frac12[L,[C,\Gamma]]$.
\par Cosequently, if $\Gamma$ is homogeneous, then $T$ is
homogeneous of degree zero.
\end{lemma}
{\bf Proof.}
We prove (e) only. By the Jacobi identity, we have\newline
$[C,[L,\Gamma]]+[L,[\Gamma,C]]-[\Gamma,[C,L]]=0$.
Consequently, $[C,[L,\Gamma]]-[L,[C,\Gamma]]+[\Gamma,L]=0$.\newline
If $\Gamma$ is homogeneous, then $[C,\Gamma]=0$ and so $[C,T]=-T$,
which means that $T$ is homogeneous of degree zero.\quad$\Box$
\begin{proposition}\quad
For any $L$-semispray $S$ on $M$, the $L$-connection
$\Gamma=[L,S]$ has\linebreak no torsion.
\end{proposition}
{\bf Proof.} 
By the Jacobi identity, we have $[L,[L,S]]-[L,[S,L]]+[S,[L,L]]=0$.
Since $[L,L]=0$, then $2[L,[L,S]]=0$; from which the result.\quad$\Box$
\begin{proposition}\quad
If $\Gamma$ is a homogeneous $L$-connection on $M$, then
$$T=\frac12[L,T^\circ]$$
\end{proposition}
{\bf Proof.} 
Taking the fact that $T$ is homogeneous of degree zero and that $[L,T]=0$
(Lemma 3), the result follows directly from Lemma 2.\quad$\Box$
\begin{corollary}\quad
For any homogeneous $L$-connection the torsion vanishes if, and
only if, its potential vanishes.
\end{corollary}
\begin{proposition}\quad
If $\Gamma$ is a homogeneous $L$-connection on $M$, then its
torsion is\linebreak given by\:
$$T=[L,[S,L]h],$$
where $S$ is the canonical $L$-spray associated with $\Gamma$.
\end{proposition}
{\bf Proof.} 
We first show that 
\begin{equation}
i_ST=[S,L]h+h
\end{equation}
Since $LS=C$ and $hS=S$, we have by (11)
$(i_ST)X=T(S,X)=v[C,hX]+v[S,LX]-L[S,hX].$
Since $[C,\Gamma]$ vanishes, then so does $[C,h]$, and so $[C,hX]=h[C,X].$
It follows that
$(i_ST)X=v[S,LX]-L[S,hX]=[S,LX]-h[S,LX]-L[S,hX]=[S,L]hX+hX$, by
(8) and (9). This proves (12). Now, using (12) together with
Proposition 3, we get 
$T=\frac12[L,[S,L]h]+\frac12[L,h]=\frac12[L,[S,L]h]+\frac12 T.$
Hence the result.\quad$\Box$
\begin{definition}\quad
Let $\Gamma$ be an $L$-connection on $M$. The strong torsion of $\Gamma$
is the L-semibasic vector $1$-form
$$t=T^\circ-\frac12[C,\Gamma].$$
\end{definition}
\begin{proposition}\quad
If $\Gamma$ is an $L$-connection on $M$, then
$$t^\circ+[C,S]-S=0,$$
where $S$ is the canonical $L$-semispray associated with $\Gamma$.
\end{proposition}
{\bf Proof.} 
We have\newline
$t^\circ=(T^\circ)^\circ-\frac12[C,\Gamma]^\circ=-\frac12[C,\Gamma]^\circ=
S-[C,S]$, by Proposition 1.\quad$\Box$
\begin{theorem}\quad
Canonical decomposition theorem\:\newline
Let $S$ be an $L$-semispray on $M$ and let $\,t\,$ be an
$L$-semibasic vector $1$-form on $M$ such that $t^\circ+[C,S]-S=0.$
There exists one and only one $L$-connection $\Gamma$ on $M$
whose canonical $L$-semispray is $S$ and whose strong torsion is~$t$~.\newline
This $L$-connection is given by
$$\Gamma=[L,S]+t$$
\end{theorem}
{\bf Proof.} 
It suffices to write the proof of Grifone [3] with the necessary
modifications. See also [1].\quad$\Box$
\begin{corollary}\quad
Let $\Gamma$ be an $L$-connection on $M$. The strong torsion of
$\Gamma$ vanishes if, and only if, $\Gamma$ is homogeneous with
no torsion.
\end{corollary}
\begin{definition}\quad
A homogeneous $L$-connection $\Gamma$ on $M$ is said to be
conservative if there exists an $L$-spray $G$ on $M$ such that 
$\;\Gamma=[L,G]$.
\end{definition}
\begin{theorem}\quad
A necessary and sufficient condition for a homogeneous
$L$-connection $\Gamma$ on $M$ to be conservative is that $\;[L,\Gamma]=0$
\end{theorem}
{\bf Proof.} 
The result follows from Theorem 3.\quad$\Box$

\Section {Curvature of an $L$-Connection}
\begin{definition}\quad
Let $\Gamma$ be an $L$-connection on $M$. The curvature of
$\Gamma$ is the vector $2$-form
$$R=-\frac12[h,h]$$
\end{definition}
\begin{lemma}\quad
The curvature $R$ of an $L$-connection $\Gamma$ on $M$ has the
following\linebreak properties\:\newline
\(a\) $R$ is $L$-semibasic.\newline
\(b\) $R=-\frac12[v,v]=\frac12[h,v]=-\frac18[\Gamma,\Gamma]$.\newline
\(c\) $R(hX,hY)=R(X,Y)$.\newline
\(d\) $R(X,Y)=-v[hX,hY]$.\newline
\(e\) $[C,R]=-\frac12[h,[C,\Gamma]]$.
\par Consequently, if $\Gamma$ is homogeneous, then $R$ is
homogeneous of degree $1$.
\end{lemma}
\begin{theorem}\quad Generalized Bianchi's identities\:\newline
The following identities relate the curvature and torsion of an
$L$-connection on $M$\:\newline
\(a\) $[L,R]=[h,T]$.\newline
\(b\) $[h,R]=0$.\newline
The last identity may be put in the form{\em:} $[[L,S],R]=-[t,R]$.
\end{theorem}
{\bf Proof.}\newline 
(a) It follows from the Jacobi identity that 
$[L,[h,h]]+[h,[h,L]]+[h,[L,h]]=0$. Hence,
$-~2[L,R]+2[h,[L,h]]=0\;$, or $\;[L,R]=[h,T]$.\newline
(b) Similarly, using the Jacobi identity, we have $[h,R]=0$.
Using the canonical decomposition theorem 3, (b) can be written
as $\; 0=[h,R]=\frac12 [\Gamma,R]=\frac12 [[L,S]+t,R].$ From which
$[[L,S],R] =-[t,R]$.\quad $\Box$
\begin{proposition}\quad
For a conservative $L$-connection $\Gamma$ on $M$, we have\:\newline 
$$R=\frac13[L,R^\circ]$$
\end{proposition}
{\bf Proof.} 
As $\Gamma$ is homogeneous, then $R$ is homogeneous of degree $1$. 
Applying Lemma~2 on $R$, we get: $R=\frac13([L,R]^\circ+[L,R^\circ]).$
But since $T=0$, then $[L,R]=0$ by Bianchi identity (a). Hence
the result.\quad$\Box$
\begin{corollary}\quad
For a conservative $L$-connection $\Gamma$ on $M$ the curvature
$R$ of $\Gamma$ vanishes if, and only if, its potential
$R^\circ$ vanishes.
\end{corollary}
\par The above result can be generalized for $L$-connections
$[L,S]$ which need not be conservative (Theorem 6 below). We
first need the following result [1].
\begin{proposition}\quad
For any $L$-semispray $S$, the curvature of the $L$-connection\linebreak
$\Gamma=[L,S]$ on $M$ is given by\:\newline
\(a\) $R=-\frac14[L,[S,h]]$.\newline
\(b\) $R=-[L,h[S,h]]$.
\end{proposition}
\begin{theorem}\quad 
The curvature $R$ of the $L$-connection $[L,S]$ vanishes if, and
only if, its potential $R^\circ$ vanishes.
\end{theorem}
{\bf Proof.} 
It is clear that if $R=0$, then $R^\circ=0$.\newline
Conversely, let $R^\circ=0$. Then, we have for all $X\in{\goth X}(M)\;$,
$0=R(S,X)=-v[S,h]X$. Hence, $v[S,h]=0\;$ or, equivalently,
$[S,h]-h[S,h]=0$. Consequently, $[L,[S,h]]-[L,h[S,h]]=0$. Using
Proposition 7, we get $-3R=0$. Hence the result.~\quad $\Box$
\begin{theorem}\quad 
Let $\Gamma$ be an $L$-connection on $M$. The following
assertions are\linebreak equivalent{\em:}\newline
(a) The horizontal distribution $z\longmapsto H_z(M)$ is
completely integrable.\newline
(b) The Nijenhuis tensor of $\Gamma$ vanishes.\newline
(c) The curvature $R$ of $\Gamma$ vanishes.
\end{theorem}
{\bf Proof.} 
The result follows from the following two equivalences:\newline
$[h,h]=0\iff R=0$\newline
$[\Gamma,\Gamma]=0\iff R=0$\quad $\Box$
\begin{definition}{\bf [7]}\quad
Let $\Gamma$ be an $L$-connection on $M$ and $R$ its curvature.
The nullity space of $R$ at a point $z\in M$ is the horizontal
subspace
$$N_R(z):=\{X\in H_z(M): i_XR=0\}$$
The map $N_R: z\longmapsto N_R(z)$ is called the nullity
distribution of $R$.
\par The nullity index of $R$ at $z\in M$ is defined by
$$\mu_R(z):=\dim N_R(z)$$
The map $\mu_R:z\longmapsto \mu_R(z)$ is called the nullity
index of $R$. Clearly, $0\le\mu_R(z)\le~n\;$ for all $z\in M$.
\end{definition}
\par We have the following result concerning the integrability
of the nullity distribution $N_R$. The proof of this result
follows the same lines as in [7] with simple necessary modifications.
\begin{theorem}\quad 
If $\mu_R$ is constant on the open set $U\subset M$, the nullity
distribution $N_R$ is completely integrable on $U$.
\end{theorem}

\Section{Almost-Complex Structure Associated}
\vspace{-0.5cm}
\begin{center}
{\Large\bf to an $L$-connection}
\end{center}
\vspace{0.4534cm}
\par In this section we show that to each $L$-connection $\Gamma$ on
$M$ there is naturally associated an almost-complex structure and we
investigate such structure.
\begin{proposition}\quad
Let $\Gamma$ be an $L$-connection on $M$. There exists a unique
vector $1$-form $F$ on $M$ such that
\begin{equation}
FL=h,\qquad Fh=-L
\end{equation}
Moreover, $F$ defines an almost-complex structure on $M$ with $LF=v$.
\end{proposition}
{\bf Proof.} 
To prove the existence of $F$ it suffices to show that it is
well-defined. Let $V\in V(M)$ and let $Y,Y'\in {\goth X}(M)$ be
such that $LY=LY'=V$. Since $Y-Y'\in Ker(L)=Im(L)$, then
$h(Y-Y')=0$. Hence, $FV=hY=hY'\,$ and the first relation of (13) is
well-defined. Similarly, the seconed relation of (13) is well-defined. 
On the other hand, $F$ is unique since it is completely
determined on both $V(M)$ and $H(M)$. 
\par Now, $F^2(LX)=F(hX)=-LX\;$ and $\,F^2(hX)=-F(LX)=-hX\;$.
Hence, $F$ is an almost-complex structure on $M$.
\par Finally, one can easily prove that $LF=v$.\quad $\Box$
\begin{definition}\quad
Let $\Gamma$ be an $L$-connection on $M$. The almost-complex
structure $F$ defined in Proposition 8 is called the
almost-complex structure associated to $\Gamma$.
\end{definition}
\begin{lemma}\quad
The almost-complex structure $F$ associated to $\Gamma$ has the
following \linebreak properties\:\newline
\(a\) $LF=v.$\newline
\(b\) $Fv=hF.$\newline
\(c\) $vF=F-Fv=F-hF=-L.$\newline
\(d\) $L[C,F]=v-\frac12[C,\Gamma].$\newline
\(e\) $h[S,h]=hF.$\newline
\end{lemma}
\begin{proposition}\quad
The almost-complex structure $F$ associated to an $L$-connection
$\Gamma$ on $M$ can be expressed by\:
$$F=h[S,h]-L$$
\end{proposition}
{\bf Proof.} 
For all $X\in{\goth X}(M)$, we have
$(h[S,h]-L)LX=h[S,h]LX=h\{[S,hLX]-h[S,LX]\}=-h^2[S,LX]=hX$, by
equation (9). Hence, $(h[S,h]-L)L=h$. Similarly, $(h[S,h]-L)h=-L$.
The result follows then from the uniqueness of $F$\quad $\Box$
\begin{proposition}\quad
The following identities relate $F$ to both $T$ and $R$\:\newline 
\(a\) $h^\star[F,F]=F\circ T+R$,\quad 
where $h^\star[F,F](X,Y)=\frac12[F,F](hX,hY)$.\newline 
\(b\) $[L,F]=i_FT-F\circ T-R$.\newline 
\(c\) $[h,F]=-i_FR-T$.
\end{proposition}
The proof of these identities is too long but simple.
So, we omit the proof.
\begin{proposition}\quad
For an $L$-connection $\Gamma$ on $M$ without torsion the
curvature $R$ of $\Gamma$ has the following forms\:\newline
\(a\) $R=h^\star[F,F]$.\newline
\(b\) $R=-[L,F]$.\newline
\(c\) $R=-[L,hF]$.\newline
\(d\) $R=-[L,h[S,h]]$.
\end{proposition}
The proof follows from Proposition 10 and Lemma 5.
\begin{theorem}\quad 
Let $\Gamma$ be an $L$-connection on $M$. The following
assertions are \linebreak equivalent\:\newline
\(a\) $F$ is completely integrable on $M$.\newline
\(b\) $T=0\;$ and $\;R=0$.\newline
\(c\) $[L,F]=0.$\newline
\(d\) $[h,F]=0.$
\end{theorem}
{\bf Proof.} 
(a)$\iff$(b): If $T=0\;$ and $\;R=0$, then $h^\star[F,F]=0$ by
Proposition 10(a). Since $[F,F](hX,hY)=F\circ h^\star[F,F](X,Y)\;$
and $\;[F,F](LX,LY)=-h^\star[F,F](X,Y)$, then $h^\star[F,F]=0$ implies
$[F,F]=0$. 
\par Coversely, if $[F,F]=0$, then by Proposition 10(a),
$0=L\circ h^\star[F,F]=v\circ T=T$. Again, using the same formula,
we also have $R=0$.\newline
(c)$\iff$(b): If $T=0\;$ and $\;R=0$, then $[L,F]=0$ by Proposition 10(b). 
\par Coversely, if $[L,F]=0$, then by Proposition 10(b),
$0=[L,F](hX,LY)=T(X,Y)$ and again, by the same formula, $R=0$.\newline
(d)$\iff$(b): If $T=0\;$ and $\;R=0$, then $[h,F]=0$ by Proposition 10(c). 
\par Coversely, if $[h,F]=0$, then by Proposition 10(c),
$0=[h,F](hX,LY)=-R(X,Y)$ and again, by the same formula, $T=0$.\quad $\Box$
\begin{proposition}\quad
Let $\Gamma$ be an $L$-connection on $M$. Any Riemannian metric
$g$ on $V_z(M)$; $z\in M$, can be extended to a Riemannian metric 
$g_\Gamma$ on $T_z(M)$ defined by\:
$$g_\Gamma(X,Y)=g(LX,LY)+g(vX,vY)\qquad\forall X,Y\in{\goth X}(M)$$
The metric $g_\Gamma$ on $M$ is completely determined by\:\newline
\(a\) $g_\Gamma(hX,vY)=0$,\newline
\(b\) $g_\Gamma(hX,hY)=g_\Gamma(LX,LY)=g(LX,LY)$.
\end{proposition}
\par The proof is straightforward.
\begin{theorem}\quad 
Let $\Gamma$ be an $L$-connection on $M$. The metric $g_\Gamma$ 
is almost-Hermitian with respect to $F$ and the K\"{a}hler form
$K_\Gamma$ of the almost-Hermitian manifold $(M,g_\Gamma,F)$ is
given by\:
$$K_\Gamma(X,Y)=g_\Gamma(X,LY)-g_\Gamma(LX,Y)\quad\forall X,Y\in
{\goth X}(M)$$
\end{theorem}
{\bf Proof.} 
One can immediately show that 
$(i_Fg_\Gamma)(X,Y)=g_\Gamma(FX,Y)+g_\Gamma(X,FY)=0.$
The K\"{a}hler form $K_\Gamma$ is:\newline
$K_\Gamma(X,Y)=g_\Gamma(FX,Y)=g_\Gamma(FhX,hY)+g_\Gamma(FhX,vY)+
g_\Gamma(FvX,hY)+g_\Gamma(FvX,vY)\newline
=-g_\Gamma(LX,hY)-g_\Gamma(LX,vY)+
g_\Gamma(hFX,hY)+g_\Gamma(hFX,vY); \hbox{ by (13) and Lemma 5}\newline
=g_\Gamma(hFX,hY)-g_\Gamma(LX,vY); \hbox{ by
Proposition 12}\newline
=g_\Gamma(LFX,LY)-g_\Gamma(LX,Y)
=g_\Gamma(vX,LY)-g_\Gamma(LX,Y)\newline
=g_\Gamma(X,LY)-g_\Gamma(LX,Y)$.
\quad $\Box$


\end{document}